\documentclass[12pt]{amsart}

\usepackage{amssymb, color}
\usepackage[dvipsnames]{xcolor}
\usepackage{mathrsfs}
\usepackage{amsmath} 
\usepackage{amsthm}
\usepackage{url} 
\usepackage{stackrel} 
\usepackage{courier}
\usepackage[all,cmtip]{xy} 
\usepackage{multicol}
\usepackage{mathtools}
\usepackage{verbatim} % \begin{comment} \end{comment}
\usepackage{lipsum}
\usepackage{hyperref}
\usepackage[normalem]{ulem}
\usepackage[makeroom]{cancel}  
\usepackage{enumitem}
\catcode`~=11 % make LaTeX treat tilde (~) like a normal character
\newcommand{\urltilde}{\kern -.15em\lower .7ex\hbox{~}\kern .04em}  
\catcode`~=13 % revert back to treating tilde (~) as an active character

\makeatletter
\newif\if@borderstar
   \def\bordermatrix{\@ifnextchar*{%
       \@borderstartrue\@bordermatrix@i}{\@borderstarfalse\@bordermatrix@i*}%
   }
   \def\@bordermatrix@i*{\@ifnextchar[{\@bordermatrix@ii}{\@bordermatrix@ii[()]}}
   \def\@bordermatrix@ii[#1]#2{%
   \begingroup
     \m@th\@tempdima8.75\p@\setbox\z@\vbox{%
       \def\cr{\crcr\noalign{\kern 2\p@\global\let\cr\endline }}%
       \ialign {$##$\hfil\kern 2\p@\kern\@tempdima & \thinspace %
       \hfil $##$\hfil && \quad\hfil $##$\hfil\crcr\omit\strut %
       \hfil\crcr\noalign{\kern -\baselineskip}#2\crcr\omit %
       \strut\cr}}%
     \setbox\tw@\vbox{\unvcopy\z@\global\setbox\@ne\lastbox}%
     \setbox\tw@\hbox{\unhbox\@ne\unskip\global\setbox\@ne\lastbox}%
     \setbox\tw@\hbox{%
       $\kern\wd\@ne\kern -\@tempdima\left\@firstoftwo#1%
         \if@borderstar\kern2pt\else\kern -\wd\@ne\fi%
       \global\setbox\@ne\vbox{\box\@ne\if@borderstar\else\kern 2\p@\fi}%
       \vcenter{\if@borderstar\else\kern -\ht\@ne\fi%
         \unvbox\z@\kern-\if@borderstar2\fi\baselineskip}%
         \if@borderstar\kern-2\@tempdima\kern2\p@\else\,\fi\right\@secondoftwo#1 $%
     }\null \;\vbox{\kern\ht\@ne\box\tw@}%
   \endgroup
   }
\makeatother

\allowdisplaybreaks 

\begin{document}

\title{$\mathfrak F$-categories and $\mathfrak F$-functors in Representation Theory II}
\author{Ben Cox}

\address{Department of Mathematics,  University of
Charleston \\
Charleston, SC 29424, USA \\
 coxbl@cofc.edu}

\maketitle

\begin{abstract}  This is a partial derivative of \cite{MR94g:17044}.  We give a list of examples/problems that some will find amusing.
\end{abstract}  
\thanks{The author is partially supported by a collaboration grant from the Simons Foundation (\#319261).  }

\maketitle 

\bibliographystyle{amsalpha} 

\setcounter{tocdepth}{0} 

%\tableofcontents 

% Check the coefficients: the roots \alpha_i  versus the poly coeff a_i's.  

\section{Introduction}\label{section:intro} 
In the tradition of I. M. Gelfand we will take some simple nontrivial examples and partially explore the consequences.  In the tradition of Grothendieck we will categorify (this appears to be a small part of what Grothendieck had in mind - who knows what he had in mind? I certainly don't.)  If you are looking for proofs, I hate to disappoint you.  As far as I can tell there are no proofs.

We use classical representation theoretic ideas found for example in  \cite{MR0450380} by J. P. Serre, in \cite{MR1153249} by W. Fulton and J. Harris, in \cite{MR2522486} by R. Goodman and N. Wallach and new ideas from C. Kassel and V. Drinfeld found in \cite{MR96e:17041}.  See also works by Bakturin, Zhelobenkov, Kirillov, Ibragimov, Lychagin, Komrakov, Vilekin, Vershik, Neretin and Vinberg.  

I view categorification as a cheap mathematical microscope and/or telescope depending on one's point of view.
\subsection{Clebsch-Gordan Coefficients and Clebsch-Gordan decomposition}
The Clebsch-Gordan decomposition for $\mathfrak{sl}(2)$ is  
\begin{equation}
F_m\otimes F_n\cong F_{m+n}\oplus \cdots \oplus F_{|m-n|}.
\end{equation}

Thus we have an $\mathfrak F$-category. 
The Clebsch-Gordan coefficients are obtained by taking a basis of $F_m$ say $\mathbf v_m, \mathbf v_{m-2}\dots ,\mathbf v_{-m}$ and equating coefficients using a non-degenerate bilinear form (see \cite{MR96e:17041} or more precisely \cite{MR2586983}). 
This isomorphism can defined on highest weight vectors by
\begin{equation}
\Phi(\mathbf v^{(m+n-2p)}_{m+n-2p})=\sum_{k=0}^{p}(-1)^{n-p}\frac{[n-p+k]![m-k]!}{[n-p]![m]!}   v^{(k-p)( 2 + m)+p^2 -k^2  + n}  \mathbf v^{(m-j)}_k\otimes
\mathbf v_{n-p+k}^{(n)}
\end{equation}

Now take your favorite finite group say $D_n$ and its mutations or avatars such as $Dic_h$.  The two have the same character theory.  But what about their categorifications?  There are a countable number of finite groups, so you have your work cut out for you.  Categorify for example results in \cite{MR513992} and \cite{MR0389024}.  Many of the groups often have geometric content. See for example \cite{MR1796706}.  Can one interpolate between categorifications?   I think the answer is yes.

\section{$\mathfrak F$-categories from other minds.}
\subsection{The category $\mathcal I$ of Enright, \cite{MR541329}}   This is the non-triangular, nonabelian but additive category 
whose indecomposable objects are 
$$
M_n\quad \text{and}\quad P_n\quad \text{ but not}\quad F_n.
$$
 Categorify $M_n$, $P_n$, functors between them and study the resulting categorical structures. 
Only part of this work has been done.  An abelian categorification of $M_n$ appears in the work of Naisse-Vaz and an additive version starting from M. Khovanov's work will appear hopefully some day.   Then the pieces will need to be put back together.

The end result should be a categorification of Enright's  Theorem in  \cite{MR541329}.
\subsection{The category $\mathcal H\mathcal T$ of \cite{MR1151617} of Howe and Tan}    We might call this non-abelian categorification or non-abelian harmonic categorification?   One needs to consult \cite{MR1151617} for background info and notation.

Consider the representation $\widetilde{(V_\lambda\otimes \bar V_\nu)}$  which has a ``basis'' $\mathbf v_n\otimes \overline{\mathbf v_k}$.  Using the action of $\mathfrak{sl}(2,\mathbb R)$ categorify this action.   Consider the module $U(\nu^+,\nu^-)$.  Categorify its structure.
\subsection{The category $\mathcal R$ of Rasskazova 
\cite{MR1302526}}
Consider the representation of Rasskazova's $V=V(\beta,\lambda,n)$, which has basis
 $$\{\mathbf v^i_j|i=1,...,n; j\in\mathbb Z\},$$
  we define the homomorphism $\varphi:\mathfrak{sl}(2,\mathbb C)\to\mathfrak{gl}(V)$, \begin{align*}
\varphi(h)(\mathbf v^i_j)=h\mathbf v^i_j&=(2j+\beta)\mathbf v^i_j\\
\varphi(e)(\mathbf v^i_j)=e\mathbf v^i_j&=\mathbf v^i_{j+1} &j\ge0,\\
\varphi(e)(\mathbf v^i_j)=e\mathbf v^i_j&=(\lambda+j\beta+j(j+1))\mathbf v^i_{j+1}+\mathbf v^{i-1}_{j+1} &j<0,\\
\varphi(f)(\mathbf v^i_j)=f\mathbf v^i_j&=-(\lambda+(j-1)\beta+j(j-1))\mathbf v^i_{j-1}-\mathbf v^{i-1}_{j-1} &j>0,\\
\varphi(f)(\mathbf v^i_j)=f\mathbf v^i_j&=-\mathbf v^i_{j-1} &j\le0.
\end{align*}

 Categorify these representations  $V(\beta,\lambda,n)$, functors between them and study the resulting categorical structures. 
I believe Rasskazova has other representations.   

Categorify any functor $F:\mathcal H\mathcal T, \mathcal R,\mathcal I\to \mathcal H\mathcal T, \mathcal R,\mathcal I$ etc.
\section{Conclusion}   There is no conclusion, but there are other partial derivatives of \cite{MR94g:17044}.  My mind reels from the possibilities.  What is the geometric content and structure? There is lots more to come into focus.  I can see parts of images right now.  I'll write about those images later.

%If you have a big head and need a big wig, i.e. a bigwig, do not fret I will likely cite you next time.  If your name ends in an "ov" then you may or may not be cited next time.  (Eg  Arkhipov or Zelmanov).
% \bibliography{math} 
% @article {MR1191610,
%    AUTHOR = {Cox, Ben},
%     TITLE = {{${\germ F}$}-categories and {${\germ F}$}-functors in the
%              representation theory of {L}ie algebras},
%   JOURNAL = {Trans. Amer. Math. Soc.},
%  FJOURNAL = {Transactions of the American Mathematical Society},
%    VOLUME = {343},
%      YEAR = {1994},
%    NUMBER = {1},
%     PAGES = {433--453},
%      ISSN = {0002-9947},
%   MRCLASS = {17B67 (17B55)},
%  MRNUMBER = {1191610},
%MRREVIEWER = {Daniel K. Nakano},
%       URL = {https://doi-org.nuncio.cofc.edu/10.2307/2154540},
%}

%\begin{thebibliography}{C18}
%
%\bibitem[Ser77]{MR0450380}
%Jean-Pierre Serre, \emph{Linear representations of finite groups},
%  Springer-Verlag, New York-Heidelberg, 1977, Translated from the second French
%  edition by Leonard L. Scott, Graduate Texts in Mathematics, Vol. 42.
% \bibitem[FH91]{MR1153249}
%William Fulton and Joe Harris, \emph{Representation theory}, Graduate Texts in
%  Mathematics, vol. 129, Springer-Verlag, New York, 1991, A first course,
%  Readings in Mathematics.
%

\end{document}